\documentclass[11pt,a4paper]{article}

\title{\bf A probabilistic proof of the fundamental gap conjecture via
the coupling by reflection}
\author{
Fuzhou Gong$^1$\footnote{Partly supported by the Key Laboratory of
RCSDS, CAS (2008DP173182), NSFC (11021161) and 973 Project (2011CB808000).}
\quad Huaiqian Li$^2$\footnote{Email: huaiqianlee@gmail.com. Partly supported by the NSFC (11401403) and the Australian Research Council (ARC) grant (DP130101302).}
\quad Dejun Luo$^1$\footnote{Email: luodj@amss.ac.cn. Partly supported by the Key Laboratory of
RCSDS, CAS (2008DP173182), NSFC (11101407) and AMSS (Y129161ZZ1).\newline
}
\vspace{3mm}\\
{\footnotesize $^1$Institute of Applied Mathematics, Academy of Mathematics and Systems Science,}\\
{\footnotesize Chinese Academy of Sciences, Beijing 100190, China}\\
{\footnotesize $^2$School of Mathematics, Sichuan University, Chengdu 610064, China}
}
\date{}

\usepackage{amssymb,amsmath,amsfonts,amsthm,array,bm,color}

\def\R{\mathbb{R}}
\def\E{\mathbb{E}}

\def\Z{\mathbb{Z}}

\def\L{\mathcal{L}}

\newcommand{\eps}{\varepsilon}
\def\ua{\uparrow}
\def\da{\downarrow}
\def\diam{\textup{diam}}

\def\d{\textup{d}}

\def\<{\langle}
\def\>{\rangle}

\def\XXint#1#2#3{{\setbox0=\hbox{$#1{#2#3}{\int}$}
\vcenter{\hbox{$#2#3$}}\kern-.5\wd0}}

\newtheorem{theorem}{Theorem}[section]
\newtheorem{lemma}[theorem]{Lemma}       
\newtheorem{corollary}[theorem]{Corollary}
\newtheorem{proposition}[theorem]{Proposition}
\newtheorem{remark}[theorem]{Remark}
\newtheorem{example}[theorem]{Example}

\begin{document}

\maketitle

\makeatletter 
\renewcommand\theequation{\thesection.\arabic{equation}}
\@addtoreset{equation}{section}
\makeatother 

\vspace{-8mm}

\begin{abstract}
Let $\Omega\subset\R^n$ be a strictly convex domain with smooth boundary and diameter $D$.
The fundamental gap conjecture claims that if $V:\bar\Omega\to\R$ is convex, then the
spectral gap of the Schr\"odinger operator $-\Delta+V$ with Dirichlet boundary condition
is greater than $\frac{3\pi^2}{D^2}$. Using analytic methods, Andrews and Clutterbuck recently
proved in [J. Amer. Math. Soc. 24 (2011), no. 3, 899--916] a more general spectral gap
comparison theorem which implies this conjecture. In the first part of the current work,
we shall give a probabilistic proof of their result via the coupling by reflection of
the diffusion processes. Moreover, we also present in the second part a simpler probabilistic
proof of the original conjecture.
\end{abstract}


\section{Introduction}

Given a bounded strictly convex domain $\Omega\subset\R^n$ with smooth boundary $\partial\Omega$, and some potential function $V:\Omega\to\R$, we consider the Schr\"odinger operator $L=-\Delta+V$ on $\Omega$ with Dirichlet boundary condition, where $\Delta$ is the standard Laplacian operator on $\R^n$. The operator $L$ has an increasing sequences of eigenvalues $\lambda_0 <\lambda_1\leq\lambda_2 \leq\cdots$, with the associated eigenfunctions $\{\phi_i\}_{i\geq0}$ which vanish on the boundary $\partial\Omega$. We assume $\phi_i$ is normalized in $L^2(\Omega,\d x)$ for all $i\geq 0$. The eigenfunction $\phi_0>0$ and eigenvalue $\lambda_0$ are also called the ground state and ground state energy, respectively.

It was conjectured by several authors (see \cite{van, Yau, Ashbaugh}) that if $V$ is a convex potential, then the difference of the first two eigenvalues of the Schr\"odinger operator $L=-\Delta+V$ (the so-called fundamental or spectral gap)  satisfies
  \begin{equation}\label{spectral-gap}
  \lambda_1-\lambda_0\geq \frac{3\pi^2}{D^2},
  \end{equation}
where $D=\diam(\Omega)$ is the diameter of $\Omega$. The readers can find in \cite[Section 1]{Andrews} a comprehensive account of the progress on various special cases of this conjecture. It was completely solved by Andrews and Clutterbuck in the recent paper \cite{Andrews}. In their proof, they introduced the notion of {\it modulus of convexity} for the potential $V$. More precisely, a function $\tilde V\in C^1([0,D/2])$ is called a modulus of convexity for $V\in C^1(\Omega)$ if for all $x,y\in\Omega, \,x\neq y$, one has
  \begin{equation}\label{modulus-convex}
  \Big\<\nabla V(x)-\nabla V(y),\frac{x-y}{|x-y|}\Big\>\geq 2\tilde V'\bigg(\frac{|x-y|}2\bigg),
  \end{equation}
where $\<\, ,\>$ and $|\cdot|$ are respectively the inner product and Euclidean norm of $\R^n$.
Intuitively, we may say that $V$ is ``more convex'' than $\tilde V$.
If the sign is reversed, then $\tilde V$ is called the {\it modulus
of concavity} for $V$. Under the condition \eqref{modulus-convex}, Andrews and
Clutterbuck proved in \cite[Theorem 1.5]{Andrews} that $\log\phi_0$ has a modulus of
concavity $\log\tilde\phi_0$: for all $x,y\in \Omega,\, x\neq y$,
  \begin{equation}\label{modulus-log-concave}
  \Big\<\nabla\log\phi_0(x)-\nabla\log\phi_0(y),\frac{x-y}{|x-y|}\Big\>
  \leq 2(\log\tilde\phi_0)'\bigg(\frac{|x-y|}2\bigg),
  \end{equation}
where $\tilde\phi_0$ is the ground state of the one dimensional Schr\"odinger operator
$\tilde L=-\frac{\d^2}{\d t^2}+\tilde V$, satisfying the Dirichlet boundary condition
on the symmetric interval $[-D/2,D/2]$. Here $\tilde V$ is extended to be an even
function on $[-D/2,D/2]$. The sharp estimate \eqref{modulus-log-concave}
enables Andrews and Clutterbuck to prove a spectral gap comparison theorem which
implies the gap conjecture \eqref{spectral-gap}; see \cite[Proposition 3.2 and Corollary 1.4]{Andrews} for details.
The proofs are slightly simplified by Y. He in the recent paper \cite{He}.

The purpose of the present paper is to give a probabilistic proof to the
fundamental gap conjecture. Here we briefly recall the literature
of estimating the spectral gap using probabilistic methods, in particular, the
coupling method which was initiated by Professors M.-F. Chen and F.-Y. Wang in the 1990s.
In \cite{ChenWang94}, the authors applied the coupling method to estimate the
first eigenvalue of the Laplacian operator on compact Riemannian manifolds
in terms of the diameter, the dimension and the lower bound on the Ricci curvature of the manifold.
They also extended their method to obtain a variational
formula for the lower bound of the spectral gap of an elliptic operator; see \cite{ChenWang97}.
The interested readers are referred to the excellent book \cite{Chen} by
M.-F. Chen, in which Chap. 2 and 3 are devoted to the estimation of the first
eigenvalue via the coupling method. A short introduction of this method can be found
in \cite[Section 6.7]{Hsu}.

This paper consists of three parts. In part one, assuming that the sharp estimate \eqref{modulus-log-concave} on the modulus of concavity of $\log\phi_0$ holds, we first provide in Subsection \ref{subsec-gap-comparison} a direct probabilistic proof of Andrews and Clutterbuck's spectral gap comparison theorem. We would like to mention that our proof is much simpler and more direct than the one given in \cite[Proposition 3.2]{Andrews}.  Then we turn to establish the log-concavity estimate \eqref{modulus-log-concave} in Subsection \ref{subsec-log-concavity}. In the arguments of these two subsections, our main tool is the coupling by reflection of diffusion processes. As an application of the estimate \eqref{modulus-log-concave}, we show in Subsection \ref{subsec-log-Sobolev} that the ground state $\phi_0$ of the Schr\"odinger operator $-\Delta+V$ satisfies the logarithmic Sobolev inequality, provided that the modulus of convexity $\tilde V$ fulfils some suitable conditions.

In part two, we restrict ourselves to the original gap conjecture \eqref{spectral-gap}, and we shall give a simpler probabilistic proof of the sharp estimate \eqref{modulus-log-concave} with $\tilde\phi_0(z)=\cos(\pi z/D)$, then the conjecture \eqref{spectral-gap} follows easily by repeating the arguments in Subsection \ref{subsec-gap-comparison}. Finally, we collect in the appendix some technical results that are needed in the paper.

\section{Probabilistic proof of Andrews and Clutterbuck's spectral gap comparison theorem}

This section is divided into three parts. We first prove in Subsection \ref{subsec-gap-comparison} that the estimate \eqref{modulus-log-concave} implies Andrews and Clutterbuck's spectral gap comparison theorem. In Subsection \ref{subsec-log-concavity}, we show that the sharp estimate \eqref{modulus-log-concave} on the modulus of log-concavity of the ground state is a consequence of the modulus of convexity \eqref{modulus-convex}. Finally, under some conditions on the modulus of convexity $\tilde V$, we prove in Subsection \ref{subsec-log-Sobolev} that the probability measure $\d\mu=\phi_0\,\d x$ on $\Omega$ satisfies the logarithmic Sobolev inequality.

\subsection{Proof of the spectral gap comparison theorem}\label{subsec-gap-comparison}

In this subsection, we start from the sharp estimate \eqref{modulus-log-concave} on the modulus of concavity of $\log\phi_0$ and give a probabilistic proof of Andrews and Clutterbuck's spectral gap comparison theorem, i.e. \cite[Theorem 1.3]{Andrews}. Compared with the arguments presented in \cite[Sections 2 and 3]{Andrews}, our proof is much simpler and avoids the technical difficulties; see e.g. the proofs of \cite[Theorem 2.1 and Proposition 3.2]{Andrews}. Indeed, we can directly prove the spectral gap comparison theorem without the intermediate result \cite[Theorem 2.1]{Andrews}.

First we introduce some notations. Throughout the paper,
$\Omega$ is a bounded strictly convex domain in $\R^n$ with smooth boundary $\partial\Omega$.
Denote by $\rho_{\partial\Omega}:\bar\Omega\to\R_+$ the distance function to the boundary
$\partial\Omega$, and $N$ the unit inward normal vector field on $\partial\Omega$.  For $r>0$,
we write $\partial_r\Omega=\{x\in \Omega:\rho_{\partial\Omega}(x)\leq r\}$ for the $r$-neighborhood
of $\partial\Omega$.
By \cite[Corollary 2.3]{Wang05}, there exists $r_0\in (0,D/2)$ such that $\rho_{\partial\Omega}$
is smooth on $\partial_{r_0}\Omega$.
Then for any $x\in \partial_{r_0}\Omega$, there exists a unique $x'\in\partial\Omega$
such that $\rho_{\partial\Omega}(x)=|x-x'|$ and $\nabla\rho_{\partial\Omega}(x)=N(x')$.
In particular, $\nabla\rho_{\partial\Omega}=N$ on the boundary $\partial\Omega$.

Now we turn to prove the spectral gap comparison theorem. We begin with the fact that $\log\phi_0$ satisfies the equation
  $$\Delta\log\phi_0+|\nabla\log\phi_0|^2=V-\lambda_0.$$
Differentiating this equation leads to
  \begin{equation}\label{comp.1}
  \Delta(\nabla\log\phi_0)+2\<\nabla\log\phi_0,\nabla(\nabla\log\phi_0)\>=\nabla V.
  \end{equation}
Inspired by equation \eqref{comp.1}, we consider the following SDE
  \begin{equation}\label{comp-SDE}
  \d X_t=\sqrt 2\,\d B_t+2\nabla\log\phi_0(X_t)\,\d t,\quad X_0=x\in\Omega.
  \end{equation}
By the properties of the ground state $\phi_0$ near the boundary $\partial\Omega$, we see that the vector field $2\nabla\log\phi_0$ fulfills the condition of Lemma \ref{sect-4-lem-1}. Hence, starting from a point $x\in\Omega$, the solution $X_t$ will not hit the boundary $\partial\Omega$ (see also \cite{Carlen, MeyerZheng}).

Next, we consider the coupling by reflection (or mirror coupling) of the process $(X_t)_{t\geq 0}$ which was first introduced by Lindvall and Rogers in \cite{Lindvall} (see \cite{ChenLi} for related studies). To this end, define the $(n\times n)$-matrix
  $$M(x,y)=I_n-2\frac{(x-y)(x-y)^\ast}{|x-y|^2},\quad x,y\in\R^n,\, x\neq y,$$
where $I_n$ is the unit matrix of order $n$ and $(x-y)^\ast$ is the transpose of the column vector $x-y$. The matrix $M(x,y)$ corresponds to the reflection mapping with respect to the hyperplane perpendicular to the vector $x-y$. For $y\in\Omega,\,y\neq x$, consider
  \begin{equation}\label{comp-coupling}
  \d Y_t=\sqrt 2\,M(X_t,Y_t)\,\d B_t+2\nabla\log\phi_0(Y_t)\,\d t,\quad Y_0=y\in\Omega.
  \end{equation}
Since $M(X_t,Y_t)$ is an orthogonal matrix, the process $Y_t$ has the same generator as $X_t$. For small $\delta>0$, we introduce the stopping times
  $$\tau_\delta=\inf\{t>0:|X_t-Y_t|=\delta\} \quad\mbox{and}\quad \sigma_\delta=\inf\{t>0:
  \rho_{\partial\Omega}(X_t)\wedge\rho_{\partial\Omega}(Y_t)=\delta\}.$$
As $\delta$ decreases to 0, $\tau_\delta$ tends to the coupling time $\tau=\inf\{t>0:X_t=Y_t\}$; we shall set as usual $Y_t=X_t$ for $t\geq\tau$. The stopping time $\sigma_\delta$ is the first time that $X_t$ or $Y_t$ reach the area
$\partial_\delta\Omega$. As the function $\log\phi_0$ is smooth with bounded derivatives on $\Omega\setminus \partial_\delta\Omega$ for any fixed $\delta>0$, we conclude that, almost surely, $\sigma_\delta<\infty$. Since the two processes $(X_t)_{t\geq 0}$ and $(Y_t)_{t\geq 0}$ do not arrive at the boundary $\partial\Omega$, it holds $\sigma_\delta\ua \infty$ almost surely as $\delta\da0$.

Define the processes
  $$\alpha_t:=\nabla\log \phi_0(X_t)-\nabla\log\phi_0(Y_t),\quad
  \beta_t:=\frac{X_t-Y_t}{|X_t-Y_t|}$$
and $F_t:=\<\alpha_t,\beta_t\>$. We mention that the process $F_t$ always makes sense, even after the coupling time $\tau$. Indeed, $F_t=0$ almost surely for $t\geq\tau$. We deduce from \eqref{comp-SDE} and \eqref{comp-coupling} that for
$t\leq \tau_\delta\wedge\sigma_\delta$,
  \begin{equation*}
  \d(X_t-Y_t)=2\sqrt2\,\beta_t\beta_t^\ast\,\d B_t+2\alpha_t\,\d t,\quad X_0-Y_0=x-y\neq0.
  \end{equation*}
Denoting by $\xi_t=|X_t-Y_t|/2$, we have
  \begin{equation}\label{comp-difference.1}
  \d\xi_t=\sqrt2\,\<\beta_t,\d B_t\>+F_t\,\d t,\quad \xi_0=|x-y|/2>0.
  \end{equation}

Suppose that $\tilde V$ is a modulus of convexity of $V:\Omega\to\R$. We extend $\tilde V$ to be an even function on $[-D/2,D/2]$. Let $\{\tilde\lambda_i\}_{i\geq0}$ (resp. $\{\tilde\phi_i\}_{i\geq0})$ be the sequence of Dirichlet eigenvalues (resp. Dirichlet eigenfunctions) of the one-dimensional Schr\"odinger operator $\tilde L=-\frac{\d^2}{\d t^2}+\tilde V$ on the interval $[-D/2,D/2]$. We need the following simple results.

\begin{lemma}\label{comparison-lem-0}
Denote by $\Psi(t)=\tilde\phi_1(t)/ \tilde\phi_0(t),\,t\in(-D/2,D/2)$. Then $\Psi$ can be smoothly extended to $[-D/2, D/2]$, and it satisfies
  $$\Psi'(D/2)=0\quad \mbox{and}\quad \Psi'(t)>0\quad  \mbox{for all } 0\leq t< D/2.$$
Moreover, there is $c_1>0$ such that $\Psi(t)\geq c_1 t$ for all $0\leq t\leq D/2$.
\end{lemma}

\begin{proof}
It is well known that $\Psi(t)=\tilde\phi_1(t)/ \tilde\phi_0(t)$ can be extended to $[-D/2,D/2]$ as a smooth function; moreover, it satisfies
  \begin{equation}\label{prep-prop-1.2}
  \Psi'' +2(\log\tilde\phi_0)'\Psi'=-(\tilde\lambda_1-\tilde\lambda_0)\Psi \quad \mbox{in } (-D/2,D/2).
  \end{equation}
From this equation we deduce that $\Psi'(D/2)=0$. Next, since $\tilde\phi_1(D/2)=0$, the fundamental theorem of calculus implies for $t\in[0,D/2)$ that
  $$-\big(\tilde\phi_0' \tilde\phi_1\big)(t)=\int_t^{\frac D2} \big(\tilde\phi_0' \tilde\phi_1\big)'(s)\,\d s
  =\int_t^{\frac D2} \big(\tilde V(s)-\tilde\lambda_0\big)\big(\tilde\phi_0 \tilde\phi_1\big)(s)\,\d s
  +\int_t^{\frac D2} \big(\tilde\phi_0' \tilde\phi_1'\big)(s)\,\d s,$$
where in the second equality we have used the eigen-equation $-\tilde\phi_0''+\tilde V \tilde\phi_0=\tilde\lambda_0\tilde\phi_0$. Similarly,
  $$-\big(\tilde\phi_0 \tilde\phi_1'\big)(t)= \int_t^{\frac D2}\big(\tilde V(s)-\tilde\lambda_1\big)\big(\tilde\phi_0 \tilde\phi_1\big)(s)\,\d s +\int_t^{\frac D2} \big(\tilde\phi_0' \tilde\phi_1'\big)(s)\,\d s.$$
The above two identities lead to
  $$\Psi'(t)=\frac{\tilde\phi_0 \tilde\phi_1'-\tilde\phi_0' \tilde\phi_1}{\tilde\phi_0^2}(t)
  =\frac{\tilde\lambda_1-\tilde\lambda_0}{\tilde\phi_0^2(t)}\int_t^{\frac D2}\big(\tilde\phi_0 \tilde\phi_1\big)(s)\,\d s>0$$
for any $0\leq t<D/2$. The last assertion is obvious.

\end{proof}

The next result is a key step to prove the spectral gap comparison theorem.

\begin{lemma}\label{comparison-lem}
Assume that the log-concavity estimate \eqref{modulus-log-concave} holds. Then
  $$\E\Psi(\xi_t)\leq e^{-(\tilde\lambda_1-\tilde\lambda_0)t}\, \Psi\bigg(\frac{|x-y|}2\bigg).$$
\end{lemma}

\begin{proof}
Since the two processes $X_t$ and $Y_t$ stay in the domain $\Omega$, it holds that $0\leq \xi_t<D/2$ almost surely. The estimate \eqref{modulus-log-concave} implies that
  \begin{equation}\label{comparison-lem.1}
  F_t\leq 2(\log\tilde\phi_0)'(\xi_t)\quad \mbox{for all } t\geq 0.
  \end{equation}
From Lemma \ref{comparison-lem-0}, we know that the function $\Psi=\tilde\phi_1/\tilde\phi_0$ has positive and bounded derivative on $[0,D/2)$. By the It\^o formula and \eqref{comp-difference.1},
  \begin{align*}
  \d\Psi(\xi_t)&=\Psi'(\xi_t)\big[\sqrt2\,\<\beta_t,\d B_t\>+F_t\,\d t\big]+\Psi''(\xi_t)\,\d t\\
  &\leq \sqrt2\,\Psi'(\xi_t)\<\beta_t,\d B_t\>+\big[2(\log\tilde\phi_0)'\, \Psi'+\Psi'' \big](\xi_t)\,\d t\\
  &=\sqrt2\,\Psi'(\xi_t)\<\beta_t,\d B_t\>-(\tilde\lambda_1-\tilde\lambda_0)\Psi(\xi_t)\,\d t,
  \end{align*}
where in the last two steps we have used \eqref{comparison-lem.1} and \eqref{prep-prop-1.2}, respectively. The above inequality is equivalent to
  $$\d\big[e^{(\tilde\lambda_1-\tilde\lambda_0)t}\Psi(\xi_t)\big]\leq \sqrt2\,e^{(\tilde\lambda_1-\tilde\lambda_0)t}\Psi'(\xi_t)\<\beta_t,\d B_t\>.$$
Integrating this inequality from 0 to $t\wedge\tau_\delta \wedge \sigma_\delta$ and taking expectation lead to
  $$\E\big[e^{(\tilde\lambda_1-\tilde\lambda_0) (t\wedge\tau_\delta \wedge \sigma_\delta)}\Psi(\xi_{t\wedge\tau_\delta \wedge \sigma_\delta})\big]\leq \Psi(\xi_0).$$
By the dominated convergence theorem, letting $\delta\to 0$ yields
  $$\E\big[e^{(\tilde\lambda_1-\tilde\lambda_0) (t\wedge\tau)}\Psi(\xi_{t\wedge\tau})\big]\leq \Psi(\xi_0),$$
since $\tau_\delta\uparrow \tau$ and $\sigma_\delta\uparrow \infty$ almost surely. Recall that $\Psi(0)=0$ and $\xi_t=0$ for all $t\geq \tau$; thus we have
  $$\E\big[e^{(\tilde\lambda_1-\tilde\lambda_0) (t\wedge\tau)}\Psi(\xi_{t\wedge\tau})\big]
  =\E\big[{\bf 1}_{\{\tau>t\}}e^{(\tilde\lambda_1-\tilde\lambda_0) t}\Psi(\xi_{t})\big]
  =\E\big[e^{(\tilde\lambda_1-\tilde\lambda_0) t}\Psi(\xi_{t})\big],$$
which leads to the desired result.
\end{proof}

We can now prove

\begin{theorem}[Spectral gap comparison]\label{2-thm-3}
The inequality \eqref{modulus-log-concave} implies $\lambda_1-\lambda_0\geq \tilde\lambda_1-\tilde\lambda_0$.
\end{theorem}

\begin{proof}
By the ground state transform,
  $$v=\frac{e^{-\lambda_1 t}\phi_1}{e^{-\lambda_0 t}\phi_0}=e^{-(\lambda_1-\lambda_0) t}\frac{\phi_1}{\phi_0}$$
is smooth on $\R_+\times \bar\Omega$ and satisfies the heat equation
  \begin{align*}
  \frac{\partial v}{\partial t}&=\Delta v+2\<\nabla\log \phi_0,\nabla v\>\quad \mbox{in } \R_+\times\Omega;\\
  v(0,\cdot) &=\phi_1/\phi_0 \quad \mbox{on } \bar\Omega.
  \end{align*}
To simplify the notation, set $v_0=\phi_1/\phi_0$, which belongs to $C^1(\bar\Omega)\cap C^\infty(\Omega)$ . By \eqref{comp-SDE} and \eqref{comp-coupling}, the function $v$ has the probabilistic representation:
  $$v(t,x)=\E v_0(X_t)\quad\mbox{and}\quad v(t,y)=\E v_0(Y_t).$$
Since $v_0$ is Lipschitz continuous on $\bar\Omega$ with a constant $K>0$, we have
  $$|v(t,x)-v(t,y)|\leq \E|v_0(X_t)-v_0(Y_t)|\leq K\E|X_t-Y_t|=2K \E\xi_t.$$
Lemma \ref{comparison-lem-0} tells us that $\Psi(z)\geq c_1 z$ for $z\in[0,D/2]$; hence
  $$|v(t,x)-v(t,y)|\leq 2Kc_1^{-1} \E\Psi(\xi_t)
  \leq 2Kc_1^{-1} e^{-(\tilde\lambda_1-\tilde\lambda_0)t}\, \Psi\bigg(\frac{|x-y|}2\bigg),$$
where the last inequality follows from Lemma \ref{comparison-lem}. Substituting the
expression of $v$ into the above inequality leads to
  $$e^{-(\lambda_1-\lambda_0)t}|v_0(x)-v_0(y)|\leq 2Kc_1^{-1} e^{-(\tilde\lambda_1-\tilde\lambda_0)t}\, \Psi\bigg(\frac{|x-y|}2\bigg),\quad \mbox{for all } x,y\in\Omega, t\geq0.$$
Since $v_0=\phi_1/\phi_0$ is not a constant,
we conclude that $\lambda_1-\lambda_0\geq \tilde\lambda_1-\tilde\lambda_0$.
\end{proof}

\subsection{Sharp log-concavity estimate of the ground state}\label{subsec-log-concavity}

In this part we shall present a probabilistic proof of Andrews and Clutterbuck's sharp estimate \eqref{modulus-log-concave} on the modulus of log-concavity of the ground state $\phi_0$. We first establish an intermediate result, namely Theorem \ref{2-thm-1}, and then deduce the desired estimate \eqref{modulus-log-concave} by applying the approximation lemma (cf. Lemma \ref{prep-lem-2}).

The next theorem is an important step to prove the sharp estimate \eqref{modulus-log-concave}.

\begin{theorem}\label{2-thm-1}
Let $\Omega\subset\R^n$ be a strictly convex domain with diameter $D$, and $\tilde V$ a modulus of convexity of $V$ on $\Omega$. Let $u_0\in C^\infty(\bar\Omega,\R_+)$ and $u:\R_+\times \bar\Omega\to\R_+$ be a smooth solution to
  \begin{equation}\label{2-thm-1.1}
  \begin{split}
  \frac{\partial u}{\partial t}&=\Delta u-Vu \quad\mbox{in }\R_+\times\Omega;\\
  u&=0 \quad\mbox{on }\R_+\times \partial\Omega \mbox{ and }
  u(0,\cdot)=u_0 \mbox{ in } \bar\Omega.
  \end{split}
  \end{equation}
Assume that $\log u_0$ admits a modulus of concavity $\psi\in C^\infty([0,D/2])$ satisfying $\psi(0)=0$ and
  \begin{equation}\label{2-thm-1.2}
  \psi''+2\psi\psi'-\tilde V'=0.
  \end{equation}
Then $\psi$ is a modulus of concavity of $\log u(t,\cdot)$ for all $t\geq0$.
\end{theorem}

To prove this theorem, we note that $\log u$ satisfies the equation
  $$\frac{\partial\log u}{\partial t}=\Delta\log u+|\nabla\log u|^2-V.$$
Differentiating this equation leads to
  \begin{equation}\label{sect-2.3}
  \frac{\partial}{\partial t}(\nabla\log u)=\Delta(\nabla\log u)
  +2\<\nabla\log u,\nabla(\nabla\log u)\>-\nabla V.
  \end{equation}
Fix $t>0$. In view of the above equation, we consider the following SDE
  \begin{equation}\label{sect-2-SDE}
  \d X_s=\sqrt 2\,\d B_s+2\nabla\log u(t-s, X_s)\,\d s,\quad 0\leq s\leq t,\quad X_0=x\in\Omega.
  \end{equation}
The main difference of this equation from \eqref{comp-SDE} is that the time parameter here takes values in the bounded interval $[0,t]$. The assertions (2) and (3) in Lemma \ref{prep-lem-1} imply that the vector field $[0,t]\times\Omega \ni (s,x)\mapsto 2\nabla\log u(s,x)$ satisfies the condition of Lemma \ref{sect-4-lem-1}. Hence, starting from any point $x\in\Omega$, the solution $X_s$ will not arrive at the boundary $\partial\Omega$. Next we consider the coupling by reflection of $(X_s)_{0\leq s\leq t}$:
  \begin{equation}\label{sect-2-coupling}
  \d Y_s=\sqrt 2\,M(X_s,Y_s)\,\d B_s+2\nabla\log u(t-s, Y_s)\,\d s,\quad 0\leq s\leq t,\quad Y_0=y,
  \end{equation}
where $y\in\Omega,\,y\neq x$. Define the stopping times
  $$\tau_\delta=\inf\{s>0:|X_s-Y_s|=\delta\} \quad\mbox{and}\quad \sigma_\delta=\inf\{s>0:
  \rho_{\partial\Omega}(X_s)\wedge\rho_{\partial\Omega}(Y_s)=\delta\},$$
for small $\delta>0$. Though the notations of the quantities are the same as those in Subsection \ref{subsec-gap-comparison}, there will be no confusion according to the context. As $\delta$ decreases to 0, $\tau_\delta$ tends to the coupling time $\tau=\inf\{s>0:X_s=Y_s\}$; we shall set $Y_s=X_s$ for $\tau\leq s\leq t$. The stopping time $\sigma_\delta$ is the first time that $X_s$ or $Y_s$ reach the area $\partial_\delta\Omega$. Since the two processes $(X_s)_{0\leq s\leq t}$ and $(Y_s)_{0\leq s\leq t}$ do not arrive at the boundary $\partial\Omega$, it holds that $\sigma_\delta\ua t$ almost surely as $\delta\da0$.

As in Subsection \ref{subsec-gap-comparison}, we define the processes
  $$\alpha_s:=\nabla\log u(t-s, X_s)-\nabla\log u(t-s, Y_s),\quad
  \beta_s:=\frac{X_s-Y_s}{|X_s-Y_s|}$$
and $F_s:=\<\alpha_s,\beta_s\>$. Remark that the process $F_s$ makes sense for all $s\in[0,t]$: indeed, $F_s=0$ almost surely for $\tau\leq s\leq t$. Furthermore,
  $$F_0=\Big\<\nabla\log u(t,x)-\nabla\log u(t,y),\frac{x-y}{|x-y|}\Big\>, \quad
  F_t=\Big\<\nabla\log u_0(X_t)-\nabla\log u_0(Y_t),\frac{X_t-Y_t}{|X_t-Y_t|}\Big\>.$$
We deduce from \eqref{sect-2-SDE} and \eqref{sect-2-coupling} that
  \begin{equation}\label{sect-2-difference}
  \d(X_s-Y_s)=2\sqrt2\,\beta_s\beta_s^\ast\,\d B_s+2\alpha_s\,\d s,\quad X_0-Y_0=x-y\neq0.
  \end{equation}

\begin{lemma}\label{sect-2-lem-3}
Assume that $\tilde V$ is a modulus of convexity of $V$. Then, for $s\leq \tau_\delta\wedge\sigma_\delta$,
  \begin{equation}\label{sect-2-lem-3.1}
  \d F_s\geq \<\beta_s,\d M_s\>+2\tilde V'(\xi_s)\, \d s,
  \end{equation}
where $\xi_s=|X_s-Y_s|/2$ and
  $$M_s=\sqrt 2\int_0^s\big[(\nabla^2\log u)(t-r, X_r)
  -(\nabla^2\log u)(t-r, Y_r)M(X_r,Y_r)\big]\,\d B_r.$$
\end{lemma}

\begin{proof}
To compute the It\^o differential of $F_s$, we shall apply the It\^o formula to $\nabla\log u(t-s, X_s)$. Remember that $t$ is fixed and $s$ is the time variable. We have
  \begin{align*}
  \d[\nabla\log u(t-s, X_s)]&=-\frac{\partial}{\partial s}(\nabla\log u)(t-s, X_s)\,\d s
  +\sqrt 2\,(\nabla^2\log u)(t-s, X_s)\,\d B_s\cr
  &\hskip13pt +\big[2\<\nabla\log u, \nabla^2\log u\>
  +\Delta(\nabla\log u)\big](t-s, X_s)\,\d s\cr
  &=\sqrt 2\,(\nabla^2\log u)(t-s, X_s)\,\d B_s+\nabla V(X_s)\,\d s,
  \end{align*}
where the last equality is due to \eqref{sect-2.3}. In the same way, for $s\leq \tau_\delta\wedge\sigma_\delta$,
  $$\d[\nabla\log u(t-s, Y_s)]=\sqrt 2\,(\nabla^2\log u)(t-s, Y_s)M(X_s,Y_s)\,\d B_s
  +\nabla V(Y_s)\,\d s.$$
Now we obtain
  \begin{equation}\label{sect-2-lem-3.2}
  \d\alpha_s=\d M_s+(\nabla V(X_s)-\nabla V(Y_s))\, \d s,\quad s\leq \tau_\delta\wedge\sigma_\delta,
  \end{equation}
where $M_s$ is a vector-valued, square integrable martingale before the stopping time $\tau_\delta\wedge\sigma_\delta$.

It remains to compute $\d\beta_s$. For $s\leq \tau_\delta\wedge\sigma_\delta$, the It\^o formula yields
  \begin{align}\label{sect-2-lem-3.3}
  \d\beta_s&=|X_s-Y_s|^{-1}\d(X_s-Y_s)+(X_s-Y_s)\,\d\big(|X_s-Y_s|^{-1}\big)\cr
  &\hskip13pt +\d(X_s-Y_s)\cdot\d\big(|X_s-Y_s|^{-1}\big).
  \end{align}
By \eqref{sect-2-difference}, we have
  \begin{align}\label{sect-2-lem-3.4}
  \d|X_s-Y_s|=\big\<\beta_s,2\sqrt2\,\beta_s\beta_s^\ast\,\d B_s\big\>+2\<\beta_s,\alpha_s\>\,\d s
  =2\sqrt2\,\<\beta_s,\d B_s\>+2F_s\,\d s,
  \end{align}
where the last equality follows from $|\beta_s|\equiv1$. Again by the It\^o formula,
  \begin{align*}
  \d\big(|X_s-Y_s|^{-1}\big)&=-\frac{\d|X_s-Y_s|}{|X_s-Y_s|^2}
  +\frac{\d|X_s-Y_s|\cdot\d|X_s-Y_s|}{|X_s-Y_s|^3}\cr
  &=-\frac{2\sqrt2\,\<\beta_s,\d B_s\>}{|X_s-Y_s|^2}-\frac{2\<\beta_s,\alpha_s\>}{|X_s-Y_s|^2}\,\d s
  +\frac{8}{|X_s-Y_s|^3}\,\d s.
  \end{align*}
Combining this identity with \eqref{sect-2-difference}, we get
  $$\d(X_s-Y_s)\cdot\d\big(|X_s-Y_s|^{-1}\big)=-\frac{8(X_s-Y_s)}{|X_s-Y_s|^3}\,\d s.$$
Substituting these computations into \eqref{sect-2-lem-3.3}, we arrive at
  \begin{equation}\label{sect-2-lem-3.5}
  \d \beta_s=\frac2{|X_s-Y_s|}(\alpha_s-\<\beta_s,\alpha_s\>\beta_s)\,\d s,
  \quad s\leq \tau_\delta\wedge\sigma_\delta.
  \end{equation}
Notice that $\beta_s$ has no martingale part.

Now by the definition of $F_s$ and \eqref{sect-2-lem-3.2}, \eqref{sect-2-lem-3.5}, we have for $s\leq \tau_\delta\wedge\sigma_\delta$,
  \begin{align*}
  \d F_s&=\<\beta_s,\d\alpha_s\>+\<\alpha_s,\d\beta_s\>+\<\d\alpha_s,\d\beta_s\>\cr
  &=\<\beta_s,\d M_s\>+\<\nabla V(X_s)-\nabla V(Y_s),\beta_s\>\,\d s
  +2\frac{|\alpha_s|^2-\<\alpha_s,\beta_s\>^2}{|X_s-Y_s|}\,\d s.
  \end{align*}
Noticing that the last term is nonnegative and $\tilde V$ is a modulus of convexity for $V$, we complete the proof.
\end{proof}

Now we are ready to prove Theorem \ref{2-thm-1}.

\begin{proof}[Proof of Theorem \ref{2-thm-1}]
Recall that $\xi_s=|X_s-Y_s|/2$. By \eqref{sect-2-lem-3.4}, it satisfies
  \begin{equation*}
  \d\xi_s=\sqrt2\,\<\beta_s,\d B_s\>+F_s\,\d s,\quad \xi_0=|x-y|/2>0.
  \end{equation*}
When $s\leq \tau_\delta\wedge\sigma_\delta$, by the It\^o formula and \eqref{2-thm-1.2},
  \begin{align}\label{2-thm-1-proof.1}
  \d\psi(\xi_s)&=\psi'(\xi_s)\big[\sqrt2\,\<\beta_s,\d B_s\>+F_s\,\d s\big]+\psi''(\xi_s)\,\d s\cr
  &=\sqrt2\,\psi'(\xi_s)\<\beta_s,\d B_s\>+\tilde V'(\xi_s)\,\d s +\psi'(\xi_s)\big[F_s-2\psi(\xi_s)\big]\,\d s.
  \end{align}
Combining Lemma \ref{sect-2-lem-3} with \eqref{2-thm-1-proof.1}, we get for $s\leq \tau_\delta\wedge\sigma_\delta$ that
  \begin{align}\label{sect-2-comparison}
  \d\big[F_s-2\psi(\xi_s)\big]&\geq \d\tilde M_s
  -2\psi'(\xi_s)\big[F_s-2\psi(\xi_s)\big]\,\d s,
  \end{align}
where $\d\tilde M_s=\<\beta_s,\d M_s\>-2\sqrt2\,\psi'(\xi_s)\<\beta_s,\d B_s\>$ is the martingale part. The above stochastic differential inequality is the key ingredient to the proof.

The inequality \eqref{sect-2-comparison} is equivalent to
  $$\d\Big(\big[F_s-2\psi(\xi_s)\big]e^{\int_0^s 2\psi'(\xi_r)\,\d r}\Big)
  \geq e^{\int_0^s 2\psi'(\xi_r)\,\d r}\d\tilde M_s.$$
Integrating from 0 to $t\wedge\tau_\delta\wedge\sigma_\delta$ leads to
  \begin{align*}
  &\big[F_{t\wedge\tau_\delta\wedge\sigma_\delta}-2\psi(\xi_{t\wedge\tau_\delta\wedge\sigma_\delta})\big]
  e^{\int_0^{t\wedge\tau_\delta\wedge\sigma_\delta} 2\psi'(\xi_r)\,\d r}\cr
  &\hskip13pt \geq\big[F_0-2\psi(\xi_0)\big]+\int_0^{t\wedge\tau_\delta\wedge\sigma_\delta}
  e^{\int_0^s 2\psi'(\xi_r)\,\d r}\d\tilde M_s.
  \end{align*}
Taking expectation and the upper limit as $\delta\to0$, we obtain
  \begin{equation}\label{sect-2-key}
  F_0-2\psi(\xi_0)\leq \limsup_{\delta\to0}\E\Big(\big[F_{t\wedge\tau_\delta\wedge\sigma_\delta}
  -2\psi(\xi_{t\wedge\tau_\delta\wedge\sigma_\delta})\big]
  e^{\int_0^{t\wedge\tau_\delta\wedge\sigma_\delta} 2\psi'(\xi_r)\,\d r}\Big).
  \end{equation}
To exchange the order of the upper limit and the expectation, by the Lebesgue--Fatou lemma (see p.17 in \cite{Yosida}), it suffices to prove that the random variables in the expectation are bounded from above. Since $\psi$ and $\psi'$ are bounded, it is enough to show the upper-boundedness of $\big(F_{t\wedge\tau_\delta\wedge\sigma_\delta}\big)_{\delta>0}$.
Indeed, by Lemma \ref{prep-lem-1}(4), we have for any $s\in[0,t]$ and $x,y\in\Omega,\, x\neq y$,
  \begin{align*}
  \Big\<\nabla\log u(s,x)-\nabla\log u(s,y),\frac{x-y}{|x-y|}\Big\>
  &=\int_0^1\Big\<\nabla^2\log u|_{(s,(1-\theta)y+\theta x)}(x-y),\frac{x-y}{|x-y|}\Big\>\,\d\theta\cr
  &\leq C_t|x-y|\leq C_tD.
  \end{align*}
Thus for all $\delta>0$, $F_{t\wedge\tau_\delta\wedge\sigma_\delta}\leq C_tD$ almost surely. Now we deduce from \eqref{sect-2-key} that
  \begin{align}\label{sect-2-key-1}
  F_0-2\psi(\xi_0)&\leq \E\Big(\limsup_{\delta\to 0}\big[F_{t\wedge\tau_\delta\wedge\sigma_\delta}
  -2\psi(\xi_{t\wedge\tau_\delta\wedge\sigma_\delta})\big]
  e^{\int_0^{t\wedge\tau_\delta\wedge\sigma_\delta} 2\psi'(\xi_r)\,\d r}\Big)\cr
  &=\E\Big(\big[F_{t\wedge\tau}-2\psi(\xi_{t\wedge\tau})\big]
  e^{\int_0^{t\wedge\tau} 2\psi'(\xi_r)\,\d r}\Big),
  \end{align}
since $\sigma_\delta\ua t$ almost surely as $\delta$ tends to 0.

Finally we show that the right hand side of \eqref{sect-2-key-1} is negative. In fact, on the set $\{\tau<t\}$,
  $$F_{t\wedge\tau}-2\psi(\xi_{t\wedge\tau})=F_\tau-2\psi(\xi_\tau)=0,$$
where the last equality is due to $F_\tau=0$ and $\psi(0)=0$; while on the set $\{\tau\geq t\}$,
  \begin{align*}
  &F_{t\wedge\tau}-2\psi(\xi_{t\wedge\tau})=F_t-2\psi(\xi_t)\cr
  &\hskip13pt =\Big\<\nabla\log u_0(X_t)-\nabla\log u_0(Y_t),\frac{X_t-Y_t}{|X_t-Y_t|}\Big\>
  -2\psi\bigg(\frac{|X_t-Y_t|}2\bigg)\leq 0,
  \end{align*}
since $\psi$ is a modulus of concavity for $\log u_0$. Therefore, the right hand side of \eqref{sect-2-key-1} is negative. As a result, $F_0\leq 2\psi(\xi_0)$, which is equivalent to
  \begin{equation*}
  \Big\<\nabla\log u(t,x)-\nabla\log u(t,y),\frac{x-y}{|x-y|}\Big\>
  \leq 2\psi\bigg(\frac{|x-y|}2\bigg).
  \end{equation*}
The proof is now complete.
\end{proof}

Having Theorem \ref{2-thm-1} in hand, we can easily obtain the sharp estimate \eqref{modulus-log-concave} on the modulus of log-concavity of $\phi_0$.

\begin{theorem}[Sharp log-concavity estimate]\label{2-thm-2}
Assume \eqref{modulus-convex}, that is, the potential $V$ has a modulus of convexity $\tilde V$. Then the log-concavity estimate \eqref{modulus-log-concave} holds.
\end{theorem}

\begin{proof}
Fix any $\kappa\in(0,1)$. Let $u_\kappa\in C^\infty(\bar\Omega)$ and $\psi_\kappa\in C^\infty([0,D/2])$ be the functions constructed in Propositions \ref{prep-prop-2} and  \ref{prep-prop-3}, respectively. Applying Theorem \ref{2-thm-1} with $u_0=u_\kappa$ and $\psi=\psi_\kappa$, we obtain that for all $x,y\in\Omega, x\neq y$,
  \begin{equation}\label{thm.0}
  \Big\<\nabla\log u(t,x)-\nabla\log u(t,y),\frac{x-y}{|x-y|}\Big\> \leq 2 \psi_\kappa\bigg(\frac{|x-y|}2\bigg)
  \end{equation}
for any $t>0$. We deduce from the approximation lemma (see Lemma \ref{prep-lem-2}) that for all $x\in\Omega$,
  $$\lim_{t\to\infty} \nabla\log u(t,x)=\lim_{t\to\infty} \frac{\nabla u(t,x)}{u(t,x)}
  =\lim_{t\to\infty} \frac{\nabla(e^{\lambda_0 t} u(t,x))}{e^{\lambda_0 t} u(t,x)} =\frac{\nabla(a_0 \phi_0(x))}{a_0 \phi_0(x)}= \nabla\log\phi_0(x).$$
Therefore, letting $t\to\infty$ in \eqref{thm.0} yields
  $$\Big\<\nabla\log \phi_0(x)-\nabla\log \phi_0(y),\frac{x-y}{|x-y|}\Big\> \leq 2 \psi_\kappa\bigg(\frac{|x-y|}2\bigg).$$
By Proposition \ref{prep-prop-3}(ii), letting $\kappa\uparrow 1$ in the above inequality yields that for all $x,y\in\Omega, x\neq y$,
  \begin{equation*}
  \Big\<\nabla\log \phi_0(x)-\nabla\log \phi_0(y),\frac{x-y}{|x-y|}\Big\> \leq 2 (\log\tilde\phi_0)'\bigg(\frac{|x-y|}2\bigg).
  \end{equation*}
The proof is complete.
\end{proof}

\subsection{Logarithmic Sobolev inequality for the ground state $\phi_0$}\label{subsec-log-Sobolev}

Now we consider the measure $\d\mu=\phi_0^2\,\d x$ on the bounded convex domain $\Omega$. Since $\phi_0$ is normalized in $L^2(\Omega,\d x)$, we have $\mu(\Omega)=1$. It is well known that $\mu$ is a symmetric measure for the second order differential
operator $A=\Delta+\nabla\log\phi_0^2\cdot\nabla$, with the domain $C_N^2(\Omega)$ of $C^2$-functions on $\Omega$ satisfying the
Neumann boundary condition.

We want to establish the logarithmic Sobolev inequality for the probability measure $\mu$ on $\Omega$. Such an inequality was first proved by L. Gross in \cite{Gross} for the standard Gaussian measure and  has been studied intensively in the past four decades. In particular, it has been shown to be equivalent to the hyper-contractivity of the corresponding diffusion semigroup. Based on the ``carr\'e du champ'' operator, Bakry and \'Emery proposed in \cite{BakryEmery} a famous criterion for the logarithmic Sobolev inequality to hold. Here is a brief introduction. Given a diffusion operator $\L$ with symmetric measure $\nu$, define
  \begin{align*}
  \Gamma(f,g)&=\frac12 \big[\L(fg)-f\L g-g\L f\big],\\
  \Gamma_2(f,g)&=\frac12 \big[\L\Gamma(f,g)-\Gamma(f,\L g)-\Gamma(g,\L f)\big],
  \end{align*}
where $f,g$ belong to some algebra $\mathcal A$ which is dense in the domain $\mathcal D(\L)$ of $\L$. Bakry and \'Emery proved in \cite[p.199, COROLLAIRE 2]{BakryEmery} that if there is a constant $C>0$ such that
  $$\Gamma_2(f,f)\geq C\Gamma(f,f),\quad \mbox{for all}\, f\in\mathcal A,$$
then the logarithmic Sobolev inequality holds:
  \begin{equation}\label{log-Sobolev-ineq}
  C\int f^2\log\frac{|f|}{\|f\|_{L^2(\nu)}}\,\d\nu\leq \int \Gamma(f,f)\,\d\nu.
  \end{equation}
Using Bakry and \'Emery's framework, we shall prove

\begin{theorem}[Logarithmic Sobolev inequality]\label{log-Sobolev}
Assume that $V\in C^1(\Omega)$ admits a modulus of convexity $\tilde V\in C^1([-D/2,D/2])$, i.e. \eqref{modulus-convex} holds. Let $\phi_0$ be the ground state of the Schr\"odinger operator $L=-\Delta+V$ such that $\d\mu=\phi_0^2\,\d x$ is a probability on $\Omega$. Assume in addition that $\tilde V$ is even, and the first Dirichlet eigenvalue $\tilde\lambda_0$ of the one-dimensional
Schr\"odinger operator $\tilde L=-\frac{\d^2}{\d t^2}+\tilde V$ on $[-D/2,D/2]$ satisfies
  \begin{equation}\label{log-Sobolev.0}
  \tilde\lambda_0>\tilde V(0).
  \end{equation}
Then $\mu$ satisfies the logarithmic Sobolev inequality with constant $C=2(\tilde\lambda_0-\tilde V(0))$.
\end{theorem}

\begin{proof}
Recall the operator $A=\Delta+\nabla\log\phi_0^2\cdot\nabla=\Delta+2\nabla\log\phi_0\cdot\nabla$. In this case, it is well known that $\Gamma(f,f)=|\nabla f|^2$ and
  $$\Gamma_2(f,f)=\|\nabla^2 f\|_{HS}^2-2\big\<(\nabla^2\log\phi_0)\nabla f,\nabla f\big\>,$$
where $\nabla^2 f$ is the Hessian of $f$ and $\|\cdot\|_{HS}$ is the Hilbert--Schmidt norm of matrices. Let $\tilde\phi_0$ be the
eigenfunction of $\tilde L$ corresponding to $\tilde\lambda_0$. Then  $\tilde\phi_0$ is strictly positive on the open interval $(-D/2,D/2)$. Since $\tilde V$ is even, it is easy to show that $\tilde\phi_0$ is also even; hence $\tilde\phi_0'(0)=0$. The condition \eqref{modulus-convex} and Theorem \ref{2-thm-2} imply that $\log\tilde\phi_0$ is a modulus of concavity for $\log\phi_0$, that is, \eqref{modulus-log-concave} holds. Let $S^{n-1}$ be the unit sphere in $\R^n$. Fix any $x\in\Omega$. Then for all $\theta\in S^{n-1}$ and $ t>0$ such that $x+t\theta\in\Omega$, we have
   $$\big\<\nabla\log\phi_0(x+t\theta)-\nabla\log\phi_0(x), \theta\big\>
  \leq 2(\log\tilde\phi_0)'(t/2).$$
Dividing both sides by $t$ and letting $t\to0$, we obtain
  \begin{equation}\label{log-Sobolev.1}
  \big\<(\nabla^2\log\phi_0)(x)\theta,\theta\big\>\leq (\log\tilde\phi_0)''(0),\quad \forall\, x\in\Omega.
  \end{equation}
On the other hand,
  $$(\log\tilde\phi_0)''(0)=\frac{\tilde\phi_0''(t)\tilde\phi_0(t)-\tilde\phi_0'(t)^2}
  {\tilde\phi_0(t)^2}\bigg|_{t=0}=\frac{\tilde\phi_0''(0)}{\tilde\phi_0(0)}.$$
Using the equation $\tilde L\tilde\phi_0=\tilde\lambda_0\tilde\phi_0$ we obtain
  $$\tilde\phi_0''(0)=(\tilde V(0)-\tilde\lambda_0)\tilde\phi_0(0).$$
Combining these results with \eqref{log-Sobolev.1}, we get $\nabla^2\log\phi_0 \leq\tilde V(0)-\tilde\lambda_0$. Now by the expressions of $\Gamma$ and $\Gamma_2$, we arrive at
  $$\Gamma_2(f,f)\geq 2(\tilde\lambda_0-\tilde V(0))\Gamma(f,f).$$
Thus the logarithmic Sobolev inequality follows from the Bakry--\'Emery criterion.
\end{proof}

\begin{remark}{\rm When $\tilde \lambda_0\leq \tilde V(0)$,  we cannot directly
apply the Bakry--\'{E}mery criterion to obtain the logarithmic Sobolev inequality. Instead,
by Lemma \ref{prep-lem-1}(4), we have $\nabla^2\log\phi_0\leq K$ for some $K\geq0$;
moreover, since $\Omega$ is bounded, it is clear that $\int_\Omega e^{\delta \rho_o^2}\,\d\mu<+\infty$
for any $\delta>0$, where $\rho_o$ is the distance function from some fixed $o\in\Omega$.
Therefore, by \cite[Theorem 1.1]{Wang01}, the logarithmic Sobolev inequality always holds
with some $C>0$. The advantage of Theorem \ref{log-Sobolev} lies in the fact that we can
get explicit constant in some special cases, which is shown in the next result.}
\end{remark}

\begin{corollary}
Assume that the potential $V\in C^1(\Omega)$ is convex. Then the measure $\d\mu=\phi_0^2\,\d x$ satisfies the logarithmic Sobolev inequality with constant $C=\frac{2\pi^2}{D^2}$.
\end{corollary}

\begin{proof}
Since $V$ is convex, its modulus of convexity is simply given by $\tilde V\equiv 0$.
The one-dimensional differential operator $-\frac{\d^2}{\d t^2}$ on the interval
$[-D/2,D/2]$ has the ground state $\tilde\phi_0(t)=\cos(\frac{\pi t}D)$
associated to the eigenvalue $\tilde\lambda_0=\frac{\pi^2}{D^2}$.
\end{proof}

Now we give a simple sufficient condition for the inequality \eqref{log-Sobolev.0} to hold.

\begin{proposition}\label{sect-2-prop-2}
Suppose that $\tilde V(0)=0$ and $\min\{\tilde V(t):t\in[-D/2,D/2]\}> -\frac{\pi^2}{D^2}$.
 Then \eqref{log-Sobolev.0} holds.
\end{proposition}

\begin{proof}
Since $\tilde V(0)=0$, it suffices to show that the eigenvalue
$\tilde\lambda_0>0$. Recall that $\tilde\lambda_0$ has the
variational expression
  $$\tilde\lambda_0=\inf\bigg\{\int_{-\frac D2}^{\frac D2}\big[|f'(t)|^2+\tilde V(t)f(t)^2\big]\,\d t:
  f\in C^1_c(-D/2,D/2) \mbox{ and } \int_{-\frac D2}^{\frac D2}f(t)^2\,\d t=1\bigg\}.$$
We fix any $ f\in C^1_c(-D/2,D/2)$ with $\int_{-D/2}^{D/2}f(t)^2\,\d t=1$.
Since $\min\{\tilde V(t):t\in[-D/2,D/2]\}>-\frac{\pi^2}{D^2}$, we can find $\delta>0$
such that $\tilde V(t)\geq -\frac{\pi^2}{D^2}+\delta$ for all $t\in[-D/2,D/2]$. Thus
  $$-\int_{-\frac D2}^{\frac D2}\tilde V(t)f(t)^2\,\d t
  \leq \bigg(\frac{\pi^2}{D^2}-\delta\bigg)\int_{-\frac D2}^{\frac D2}f(t)^2\,\d t
  =\frac{\pi^2}{D^2}-\delta \leq \int_{-\frac D2}^{\frac D2}|f'(t)|^2\,\d t-\delta,$$
where the last inequality follows from the fact that $\frac{\pi^2}{D^2}$ is the
first Dirichlet eigenvalue of $-\frac{\d^2}{\d t^2}$ on $[-D/2,D/2]$. As a result,
  $$\int_{-\frac D2}^{\frac D2}\big[|f'(t)|^2+\tilde V(t)f(t)^2\big]\,\d t\geq \delta,$$
which implies that $\tilde\lambda_0\geq\delta>0$.
\end{proof}

In the following, we shall give an example where the potential $V$ is not convex.

\begin{example}
{\rm Let $\beta>0$ be a constant. Consider the double-well potential $\tilde V(t)=-\frac12 t^2
+\beta^2 t^4$ on the interval $[-1/(\sqrt 2 \beta),1/(\sqrt 2 \beta)]$. It is clear that
$0\geq \tilde V(t)\geq -1/(16\beta^2)$ for all $|t|\leq 1/(\sqrt 2 \beta)$, and the minimum
is attained at $t=\pm 1/(2\beta)$. Note that $D=\sqrt 2/\beta$ in this case.
Hence, by Proposition \ref{sect-2-prop-2}, the eigenvalue $\tilde\lambda_0$ of the operator
$\tilde L=-\frac{\d^2}{\d t^2}+\tilde V$ is positive when $\beta>(8\pi^2)^{-1/4}$.

As in \cite[Section 5]{Andrews}, we now define
  $$V(x)=\tilde V(|x|)+c\sum_{i=2}^n x_i^2,\quad x\in\R^n \mbox{ and }|x|\leq \frac1{\sqrt 2\, \beta}.$$
When $c$ is large enough, $V$ is a double-well potential which coincides with $\tilde V$
on the $x_1$-axis. It can be checked that $\tilde V$ is a modulus of convexity for $V$.}
\end{example}

Finally, we want to mention that, using the approximation argument developed in \cite{GongLiuLuo},
we can extend Theorem \ref{log-Sobolev} to the case of the whole $\R^n$ and the abstract Wiener space with the Ornstein--Uhlenbeck
operator. However, we do not want to go into details here since it deviates from the main point of the current paper.

\section{A more direct probabilistic proof of the fundamental gap conjecture}

In this section we restrict ourselves to the original gap conjecture. We assume
the potential function $V:\Omega\to\R$ is convex, and then present a direct
proof of the log-concavity estimate of the ground state $\phi_0$. Based on this estimate, we can give a shorter proof of the fundamental gap conjecture \eqref{spectral-gap} by following the arguments in Subsection \ref{subsec-gap-comparison}.

We intend to prove that if $V$ is convex, then the ground state $\phi_0$
satisfies
  \begin{equation}\label{sect-3.1}
  \Big\<\nabla\log\phi_0(x)-\nabla\log\phi_0(y),\frac{x-y}{|x-y|}\Big\>\leq
  -\frac{2\pi}D \tan\bigg(\frac{\pi|x-y|}{2D}\bigg),\quad \forall\,x,y\in\Omega, x\neq y.
  \end{equation}
We still use the notations introduced in Subsection \ref{subsec-gap-comparison}. In particular, we have
  $$F_0=\Big\<\nabla\log\phi_0(x)-\nabla\log \phi_0(y),\frac{x-y}{|x-y|}\Big\>.$$

\begin{lemma}\label{sect-3-lem-1}
Assume that the potential $V:\bar\Omega\to\R$ is convex. Then for $t\leq\tau_\delta \wedge \sigma_\delta$,
  \begin{equation}\label{sect-3-lem-1.1}
  \d F_t\geq \<\beta_t,\d M_t\>,
  \end{equation}
where
  $$M_t=\sqrt 2\int_0^t\big[(\nabla^2\log\phi_0)(X_s)
  -(\nabla^2\log\phi_0)(Y_s)M(X_s,Y_s)\big]\,\d B_s.$$
\end{lemma}

\begin{proof}
The proof is similar to that of Lemma \ref{sect-2-lem-3}, and hence we omit it to save space.
\end{proof}

Let $\tilde\phi_{D,0}(z)=\cos\frac{\pi z} D,\,z\in[-D/2,D/2]$ be the first Dirichlet eigenfunction of
the operator $-\frac{\d^2}{\d z^2}$ on the interval $[-D/2,D/2]$. Here and below we
write $\tilde\phi_{D,0}$ instead of $\tilde\phi_{0}$ to stress the dependence
on the length of the interval $[-D/2,D/2]$. For simplification of notations,
set $\psi_D(z)=(\log\tilde\phi_{D,0})'(z)=-\frac{\pi}D \tan \frac{\pi z}{D}$, which is
well defined on $(-D/2,D/2)$. Note that $\psi_D$ explodes at $z=\pm D/2$.
Thus we first take $D_1>D$ and consider $\tilde\phi_{D_1,0}$ and $\psi_{D_1}$.
Since $\psi_{D_1}$ is smooth on $[0,D/2]$ with bounded derivatives, it satisfies
  \begin{equation}\label{sect-3.3}
  \psi''_{D_1}+2\psi_{D_1}\psi'_{D_1}=0.
  \end{equation}
Now we are ready to prove

\begin{theorem}[Modulus of log-concavity]\label{sect-3-thm-1}
Assume that the potential function $V:\Omega\to\R$ is convex. Then for all $x,y\in\Omega$
with $x\neq y$, the estimate \eqref{sect-3.1} holds.
\end{theorem}

\begin{proof}
We follow the idea of the proof of Theorem \ref{2-thm-1}. Fix $\delta>0$ small enough and $D_1>D$.
When $t\leq \tau_\delta \wedge \sigma_\delta$, by \eqref{comp-difference.1} and the It\^o formula,
  \begin{align*}
  \d\psi_{D_1}(\xi_t)&=\psi'_{D_1}(\xi_t)\big[\sqrt2\,\<\beta_t,\d B_t\>+F_t\,\d t\big]+\psi''_{D_1}(\xi_t)\,\d t\\
  &=\sqrt2\,\psi'_{D_1}(\xi_t)\<\beta_t,\d B_t\>
  +\psi'_{D_1}(\xi_t)\big[F_t-2\psi_{D_1}(\xi_t)\big]\,\d t,
  \end{align*}
where the second equality follows from \eqref{sect-3.3}.
Combining Lemma \ref{sect-3-lem-1} and the above identity, we get for $t\leq \tau_\delta
\wedge \sigma_\delta$,
  \begin{align}\label{sect-3-thm-1.1}
  \d\big[F_t-2\psi_{D_1}(\xi_t)\big]&\geq \d\tilde M_t
  -2\psi'_{D_1}(\xi_t)\big[F_t-2\psi_{D_1}(\xi_t)\big]\,\d t,
  \end{align}
in which $\d\tilde M_t=\<\beta_t,\d M_t\>-2\sqrt2\,\psi'_{D_1}(\xi_t)\<\beta_t,\d B_t\>$
is the martingale part.

The inequality \eqref{sect-3-thm-1.1} is equivalent to
  $$\d\Big(\big[F_t-2\psi_{D_1}(\xi_t)\big]e^{\int_0^t 2\psi'_{D_1}(\xi_s)\,\d s}\Big)
  \geq e^{\int_0^t 2\psi'_{D_1}(\xi_s)\,\d s}\d\tilde M_t.$$
Integrating from 0 to $t\wedge \tau_\delta \wedge\sigma_\delta$ leads to
  \begin{align*}
  &\big[F_{t\wedge\tau_\delta \wedge\sigma_\delta}-2\psi_{D_1}(\xi_{t\wedge\tau_\delta \wedge\sigma_\delta})\big]
  e^{\int_0^{t\wedge\tau_\delta \wedge\sigma_\delta} 2\psi'_{D_1}(\xi_s)\,\d s}\cr
  &\hskip6pt \geq\big[F_0-2\psi_{D_1}(\xi_0)\big]+\int_0^{t\wedge\tau_\delta \wedge\sigma_\delta}
  e^{\int_0^s 2\psi'_{D_1}(\xi_r)\,\d r}\d\tilde M_s.
  \end{align*}
Taking expectation on both sides, we obtain
  \begin{align*}
  F_0-2\psi_{D_1}(\xi_0)&\leq \E\Big(\big[F_{t\wedge\tau_\delta \wedge\sigma_\delta}
  -2\psi_{D_1}(\xi_{t\wedge\tau_\delta \wedge\sigma_\delta})\big]
  e^{\int_0^{t\wedge\tau_\delta \wedge\sigma_\delta} 2\psi'_{D_1}(\xi_s)\,\d s}\Big).
  \end{align*}
By Brascamp and Lieb's result (see \cite[Theorem 6.1]{Brascamp}), the ground
state $\phi_0$ is log-concave, which implies the random variables $F_{t\wedge\tau_\delta
\wedge\sigma_\delta}$ are non-positive almost surely. Therefore,
  \begin{align}\label{sect-3-thm-1.2}
  F_0-2\psi_{D_1}(\xi_0) \leq -2\,\E\Big(\psi_{D_1}(\xi_{t\wedge\tau_\delta \wedge\sigma_\delta})
  e^{\int_0^{t\wedge\tau_\delta \wedge\sigma_\delta} 2\psi'_{D_1}(\xi_s)\,\d s}\Big).
  \end{align}
Moreover, by \cite[Example 5]{Lindvall}, the log-concavity of $\phi_0$ implies
the coupling of processes $(X_t)_{t\geq 0}$ and $(Y_t)_{t\geq 0}$ is successful,
that is, $\tau<+\infty$ almost surely.
Since $\psi'_{D_1}(z)=-\frac{\pi^2}{D_1^2}\sec^2(\frac{\pi z}{D_1})$ is negative for
$z\in[0,D/2]$, the term $e^{\int_0^{t\wedge\tau_\delta \wedge\sigma_\delta} 2\psi'_{D_1}(\xi_s)\,\d s}
\leq 1$ for all $t>0$. Note that $\psi_{D_1}$ is a bounded function on $[0,D/2]$. By the
dominated convergence theorem, letting $t\to+\infty$ and $\delta\to 0$ in
\eqref{sect-3-thm-1.2} gives us
  $$F_0-2\psi_{D_1}(\xi_0) \leq -2\,\E\Big(\psi_{D_1}(\xi_{\tau})
  e^{\int_0^{\tau} 2\psi'_{D_1}(\xi_s)\,\d s}\Big)=0.$$
Thus we obtain \eqref{sect-3.1} with $D$ being replaced by $D_1$. Letting $D_1$ tend to $D$ yields the desired inequality.
\end{proof}

\begin{remark}{\rm
Here we briefly explain why we cannot use this method to deal with the general case,
i.e., the modulus of convexity \eqref{modulus-convex} implies the log-concavity estimate
\eqref{modulus-log-concave} of the ground state. Indeed, in the above proof,
we make use of the fact that $\psi'_{D_1}(z)$ is non-positive, which in turn implies the
random variables $e^{\int_0^{t\wedge\tau_\delta \wedge\sigma_\delta} 2\psi'_{D_1}(\xi_s)\,\d s}$
are bounded, hence uniformly integrable. However, it is not clear whether
 such a result still holds in the general case.}
\end{remark}

With the log-concavity estimate \eqref{sect-3.1} on the ground state $\phi_0$, we can present a probabilistic proof of the fundamental gap conjecture \eqref{spectral-gap}, following the arguments in Subsection \ref{subsec-gap-comparison}. Indeed, it is a special case of Theorem \ref{2-thm-3}; hence we omit it here and only mention that the function $\Psi(z)$ in Lemma \ref{comparison-lem} is replaced by $\sin(\pi z/D),\, 0\leq z\leq D/2$.

\appendix

\section{Some technical results}

Let $\Omega$ be a bounded smooth convex domain. Given a time-dependent smooth vector field $b:\R_+\times\Omega\to\R^n$ and a
standard Brownian motion $(B_t)_{t\geq0}$ on $\R^n$, we consider the SDE
  \begin{equation}\label{sect-4.1}
  \d X_t=\sqrt 2\, \d B_t+b(t,X_t)\,\d t,\quad X_0=x\in\Omega.
  \end{equation}
We first give a sufficient condition to ensure that the process $(X_t)_{t\geq 0}$ stays in the domain $\Omega$.

\begin{lemma}\label{sect-4-lem-1}
Assume that the vector field $b$ satisfies
  \begin{equation}\label{sect-4-lem-1.1}
  \liminf_{\rho_{\partial\Omega}(x)\to 0}\inf_{t\geq0}\rho_{\partial\Omega}(x)
  \<b(t,x),\nabla\rho_{\partial\Omega}(x)\> >1.
  \end{equation}
Then for any $x\in\Omega$, almost surely, $X_t\in\Omega$ for all $t\geq0$.
\end{lemma}

\begin{proof}
First we choose a smooth function $\rho:\bar\Omega\to\R$ such that $\rho(x)>0$ for all $x\in\Omega$, and $\rho(x)= \rho_{\partial\Omega}(x)$ for all $x\in\partial_{r_0} \Omega$ (see the beginning of Subsection \ref{subsec-gap-comparison} for its definition). Here we may assume that $r_0>0$ is small enough such that $\rho_{\partial\Omega}$ is smooth on $\partial_{2r_0}\Omega$.  $\rho$ can be chosen as $f\circ\rho_{\partial\Omega}$, where $f:\R_+ \to \R_+$ is a smooth increasing function such that $f(t)=t$ for $t\in[0,r_0]$ and $f(t)=3r_0/2$ for $t\geq 2r_0$.  Then there is $c_0>0$ such that
  \begin{equation}\label{sect-4-lem-1.2}
  \Delta\rho(x)\geq -c_0,\quad \mbox{for all } x\in\bar\Omega.
  \end{equation}
By the It\^o formula,
  $$\d\rho(X_t)=\sqrt2\,\<\nabla\rho(X_t),\d B_t\>+\<\nabla\rho(X_t),b(t, X_t)\>\,\d t
  +\Delta\rho(X_t)\,\d t.$$
It is enough to study the behavior of $X_t$ near the boundary $\partial\Omega$.
When $X_t\in\partial_{r_0}\Omega$, by \eqref{sect-4-lem-1.2}, we have (cf. \cite[(2.2)]{Wang09})
  $$\d\rho_{\partial\Omega}(X_t)\geq \sqrt2\,\d W_t+\<\nabla\rho_{\partial\Omega}(X_t),
  b(t, X_t)\>\,\d t-c_0\,\d t,$$
where $W_t$ is a one-dimensional Brownian motion. By \eqref{sect-4-lem-1.1},
we can find $r_1\in(0,r_0]$ such that for all $(t,x)\in \R_+\times\partial_{r_1}\Omega$,
  $$\<\nabla\rho_{\partial\Omega}(x),b(t, x)\> -c_0
  \geq \frac1{\rho_{\partial\Omega}(x)}.$$
Thus, if $X_t\in\partial_{r_1}\Omega$, then we have
  $$\d\rho_{\partial\Omega}(X_t)\geq \sqrt2\,\d W_t
  +\frac1{\rho_{\partial\Omega}(X_t)}\,\d t.$$
Now we can apply \cite[Chap. VI, Theorem 3.1]{IkedaWatanabe} to conclude that, almost surely,
$\rho_{\partial\Omega}(X_t)>0$ for all $t\geq 0$.
\end{proof}

The next result is concerned with the well known properties of solutions to heat equations of Schr\"odinger operators.

\begin{lemma}\label{prep-lem-1}
Let $u_0\in C^\infty(\bar\Omega)$ be positive in $\Omega$ such that $u_0=0$ and $\nabla u_0\neq0$ on $\partial\Omega$.
Let $u:\R_+\times \bar\Omega\to\R_+$ be a smooth solution to
  \begin{equation}\label{prep-lem-1.1}
  \begin{split}
  \frac{\partial u}{\partial t}&=\Delta u-Vu \quad\mbox{in }\R_+\times\Omega;\\
  u&=0 \quad\mbox{on }\R_+\times \partial\Omega \mbox{ and }
  u(0,\cdot)=u_0 \mbox{ in } \bar\Omega.
  \end{split}
  \end{equation}
Then the solution $u$ verifies that
\begin{itemize}
\item[\rm(1)] for any $T>0$, $\theta_T:=\inf_{[0,T]\times\partial\Omega}|\nabla u|>0$;
\item[\rm(2)] for every $x\in\partial\Omega$, $\nabla u(t,x)=|\nabla u(t,x)|N(x)$;
\item[\rm(3)] $\lim_{\rho_{\partial\Omega}(x)\to0}\frac{u(t,x)}
{|\nabla u(t,x)|\rho_{\partial\Omega}(x)}=1$ uniformly in $t\in[0,T]$;
\item[\rm(4)] for any $T>0$, there exists $C_T\geq 0$ such that $\nabla^2\log u|_{(t,x)}(y,y)\leq
C_T|y|^2$ for all $t\in[0,T],\,x\in\Omega$ and $y\in\R^n$. Here $\nabla^2\log u$ is the Hessian
matrix of $\log u$.
\end{itemize}
\end{lemma}

\begin{proof}
The assertions (1) and (2) are known; see the beginning of the proof of \cite[Lemma 4.2]{Andrews}.
(3) is a consequence of (2). The proof of the last assertion is a little technical (cf. \cite[Lemma 4.2]{Andrews}); we omit it here to save space.
\end{proof}

We also need the following important approximation lemma (see \cite[Lemma 2.1]{Lee} or \cite[Lemma 3]{He} for its proof).

\begin{lemma}[Approximation lemma]\label{prep-lem-2}
Let $u$ be the same as in Lemma \ref{prep-lem-1}. Define $a_0=\int_\Omega u_0(x)\phi_0(x)\,\d x$, where $\phi_0$ is the ground state of $-\Delta+V$ on $\Omega$. Then there exists a constant $C>0$ such that
  \begin{equation}\label{prep-lem-2.1}
  \big\|e^{\lambda_0 t} u(t,\cdot)-a_0\phi_0\big\|_{C^k(\Omega)}\leq Ce^{-(\lambda_1-\lambda_0)t},\quad k\in \Z_+.
  \end{equation}
Here $\|\cdot\|_{C^0(\Omega)}$ is the supremum norm of functions on $\Omega$.
\end{lemma}

Recall that $\{\tilde\lambda_i\}_{i\geq0}$ (resp. $\{\tilde\phi_i\}_{i\geq0})$ is the sequence of eigenvalues (resp. eigenfunctions) of the one-dimensional Schr\"odinger operator $\tilde L=-\frac{\d^2}{\d t^2}+\tilde V$ on the interval $[-D/2,D/2]$ with the Dirichlet boundary condition.

\begin{proposition}\label{prep-prop-1}
There exists a positive constant $c_1$ such that
  \begin{equation}\label{prep-prop-1.1}
  -c_1-\frac{2}{D-2t}\leq (\log\tilde\phi_0)'(t)\leq c_1-\frac{2}{D-2t},\quad 0\leq t< D/2.
  \end{equation}
\end{proposition}

\begin{proof}
We follow the proof of \cite[Proposition A1]{He}. It is clear that $\tilde\phi_0'(D/2)<0$. Therefore, the function
  $$f(t):=\frac{\tilde\phi_0(t)}{(D/2) -t},\quad t\in[0,D/2)$$
can be smoothly extended to $[0,D/2]$, and satisfies that $f(t)>0$ for all $t\in [0,D/2]$. As a result,
  $$(\log\tilde\phi_0)'(t)=\frac{\tilde\phi_0'(t)}{\tilde\phi_0(t)}=\frac{f'(t)}{f(t)}-\frac1{(D/2)-t},\quad t\in [0,D/2).$$
Letting $c_1:=\sup_{0\leq t\leq D/2} |f'(t)/f(t)|$, which is finite, we obtain the first result.
\end{proof}

Next, we follow the arguments in \cite[Appendix B]{He} (see also \cite{Lu}) to show that, for any $\kappa<1$ in a small neighborhood of $1$, there exists a function $u_\kappa\in C^\infty(\bar\Omega, \R_+)$ with non-zero gradient on $\partial\Omega$, which admits $\kappa \log \tilde\phi_0$ as its modulus of log-concavity. This will be done in several steps. Recall that $\rho\in C^\infty(\bar\Omega,\R_+)$ satisfies $\rho|_{\partial_{r_0}\Omega}= \rho_{\partial\Omega}|_{\partial_{r_0}\Omega}$.

\begin{lemma}\label{prep-lem-4}
For any $\theta_0\in (0,1)$, we can find $\eps_0\in (0, r_0]$ such that
  $$\Big\<\nabla\rho(x),\frac{y-x}{|y-x|}\Big\> \geq \theta_0,\quad \mbox{for all } x,y \in \Omega \mbox{ with } |x-y|\geq D- \eps_0.$$
\end{lemma}

\begin{proof}
We consider the function
  $$F(x,y)=\Big\<\nabla \rho(x),\frac{y-x}{|y-x|}\Big\>,\quad x,y\in \bar\Omega \mbox{ and } x\neq y,$$
which is continuous on the closed set $\{(x,y)\in \bar\Omega\times \bar\Omega: |x-y|\geq D/2\}$. Let $\bar S$ be the closed subset of $\bar\Omega\times \bar\Omega$ consisting of pairs of points such that their distance is exactly $D$, i.e. $\bar S=\{(x,y)\in \bar\Omega\times \bar\Omega: |x-y|=D\}$. If $(x_0,y_0)\in \bar S$, then we have $x_0,y_0\in \partial\Omega$ and $\nabla\rho(x_0)=\nabla \rho_{\partial\Omega}(x_0) =N(x_0)=\frac{y_0-x_0}{|y_0-x_0|}$, thus $F(x_0,y_0)=1$.

For sufficiently small $\eps>0$, we define
  $$F_\eps=\min\big\{F(x,y):x,y\in \bar\Omega \mbox{ with } |x-y| \geq D-\eps\big\}.$$
Then $F_\eps$ increases as $\eps$ tends to 0. Note that if $x\in \partial\Omega$, then $|\nabla\rho(x)|=|N(x)| =1$, which implies $F_\eps\leq 1$. Thus the limit $\bar\theta:=\lim_{\eps\to 0} F_\eps$ exists and $\bar\theta \leq 1$. If $\bar\theta <1$, then we can find a sequence $\{(x_n,y_n)\}_{n\geq 1}\subset \bar\Omega\times \bar\Omega$ such that $|x_n-y_n|\to D$ and $F(x_n,y_n)\leq \bar\theta$. Extracting a subsequence if necessary, we may assume $(x_n,y_n)\to (\bar x,\bar y)$ as $n\to \infty$, hence $(\bar x,\bar y)\in \bar S$ and $F(\bar x,\bar y)\leq \bar\theta <1$. This is a contradiction. Therefore, $\lim_{\eps\to 0} F_\eps =1$ which yields the desired result.
\end{proof}

\begin{lemma}\label{prep-lem-5}
For any $\kappa\in (0,1)$, there exists a positive constant $c_2$, depending on $\kappa$, such that for all $x,y\in\Omega, x\neq y$,
  \begin{equation}\label{lem-2.0}
  \Big\<\nabla\log\rho(y)- \nabla\log\rho(x),\frac{y-x}{|y-x|}\Big\>\leq 2\kappa (\log\tilde\phi_0)'\bigg(\frac{|x-y|}2\bigg) +c_2|x-y|.
  \end{equation}
\end{lemma}

\begin{proof}
We fix an arbitrary $\kappa\in (0,1)$. Take $\theta_0\in (\kappa,1)$ and $\alpha>0$ such that $\kappa=4\theta_0/(4+\alpha)$. Let $\eps_0>0$ be determined as in Lemma \ref{prep-lem-4} for this $\theta_0$. Set $\eps_1=\min\{\eps_0,\alpha/2c_1\}$, where $c_1$ is given in Proposition \ref{prep-prop-1}.

Fix any $x,y\in\Omega, x\neq y$. There are two different cases. First, if $|y-x|>D-\eps_1$, then we can apply Lemma \ref{prep-lem-4} to get
  \begin{align}\label{lem-2.0.5}
  \Big\<\nabla\log\rho(y)- \nabla\log\rho(x),\frac{y-x}{|y-x|}\Big\>
  &=-\bigg[\Big\<\frac{\nabla\rho(y)}{\rho(y)},\frac{x-y}{|x-y|}\Big\>
  +\Big\<\frac{\nabla\rho(x)}{\rho(x)},\frac{y-x}{|y-x|}\Big\>\bigg]\cr
  &\leq - \theta_0\bigg[\frac1{\rho(y)}+ \frac1{\rho(x)}\bigg].
  \end{align}
Suppose the straight line passing from $x$ to $y$ intersects the boundary $\partial\Omega$ first at the point $\hat x$ and then at $\hat y$. We have
  \begin{align*}
  \rho_{\partial\Omega}(x)+|x-y|+\rho_{\partial\Omega}(y) &\leq |\hat x-x|+|x-y|+|y-\hat y|=|\hat x-\hat y|\leq D.
  \end{align*}
Hence $\rho_{\partial\Omega}(x)+\rho_{\partial\Omega}(y)\leq D-|x-y|<\eps_1$, thus $x,y\in \partial_{\eps_1}\Omega \subset \partial_{r_0}\Omega$, which implies $\rho(x)=\rho_{\partial\Omega}(x)$ and $\rho(y)=\rho_{\partial\Omega}(y)$. Therefore,
  \begin{equation}\label{lem-2.1}
  \rho(x)+\rho(y)\leq D-|x-y|< \eps_1.
  \end{equation}
As a result,
  $$(\rho(x)+\rho(y))(D-|x-y|)\geq (\rho(x)+\rho(y))^2\geq 4\rho(x)\rho(y),$$
which, together with Proposition \ref{prep-prop-1}, implies
  \begin{align*}
  -\bigg[\frac1{\rho(y)}+ \frac1{\rho(x)}\bigg] &\leq -\frac4{D-|x-y|} =-\frac{4}{4+\alpha}\bigg(\frac4{D-|x-y|} +\frac\alpha{D-|x-y|}\bigg)\\
  &\leq \frac{4}{4+\alpha}\bigg[2(\log\tilde\phi_0)'\bigg(\frac{|x-y|}2\bigg)+2c_1-\frac\alpha{D-|x-y|}\bigg]\\
  &\leq \frac{8}{4+\alpha}(\log\tilde\phi_0)'\bigg(\frac{|x-y|}2\bigg),
  \end{align*}
where the last inequality follows from \eqref{lem-2.1} and $\eps_1\leq \alpha/2c_1$. Combining this inequality with \eqref{lem-2.0.5}, we arrive at
  \begin{align*}
  \Big\<\nabla\log\rho(y)- \nabla\log\rho(x),\frac{y-x}{|y-x|}\Big\>
  &\leq 2\kappa(\log\tilde\phi_0)'\bigg(\frac{|x-y|}2\bigg).
  \end{align*}
Therefore, we obtain \eqref{lem-2.0} with $c_2=0$ in this case.

Next, we consider the second case that $0<|x-y|\leq D-\eps_1$. Noticing that
  $$(\log\tilde\phi_0)''(t)=\frac{\tilde\phi_0 \tilde\phi_0''-(\tilde\phi_0')^2}{\tilde\phi_0^2}(t)
  =\tilde V(t)-\tilde\lambda_0 -\bigg(\frac{\tilde\phi_0'}{\tilde\phi_0}\bigg)^2(t),$$
we can find $\bar c>0$ such that
  \begin{equation*}
  \big|(\log\tilde\phi_0)''(t)\big|\leq \bar c,\quad \mbox{for all } t\in[0,(D-\eps_1)/2].
  \end{equation*}
Moreover, $\tilde\phi_0'(0)=0$ since $\tilde\phi_0$ is an even function. Thus, for any $t\in[0,(D-\eps_1)/2]$,
  \begin{equation}\label{lem-2.2}
  \big|(\log\tilde\phi_0)'(t)\big|=\big|(\log\tilde\phi_0)'(t)-(\log\tilde\phi_0)'(0)\big| \leq \bar c t.
  \end{equation}
Finally, it is clear (cf. Lemma \ref{prep-lem-1}(4)) that $\nabla^2\log\rho\leq \hat c$ for some $\hat c\geq 0$; hence
  \begin{align*}
  \Big\<\nabla\log\rho(y)- \nabla\log\rho(x),\frac{y-x}{|y-x|}\Big\>
  &=\frac{1}{|x-y|}\int_0^1 \big\<\nabla^2\log\rho|_{x+t(y-x)}(y-x),y-x\big\>\,\d t\\
  &\leq \hat c|x-y|.
  \end{align*}
Combining the above inequality with \eqref{lem-2.2}, we get \eqref{lem-2.0} with $c_2=\bar c+\hat c$ in the second case.
\end{proof}

We can now prove

\begin{proposition}\label{prep-prop-2}
Fix any $\kappa\in(0,1)$ and let $c_2$ be given in Lemma \ref{prep-lem-5}. Define the function
  $$u_0(x)=e^{-c_2|x|^2/2}\rho(x),\quad x\in \bar\Omega.$$
Then for all $x,y\in\Omega, x\neq y$, we have
  $$\Big\<\nabla\log u_0(y)- \nabla\log u_0(x),\frac{y-x}{|y-x|}\Big\>\leq 2\kappa (\log\tilde\phi_0)'\bigg(\frac{|x-y|}2\bigg).$$
\end{proposition}

\begin{proof}
By the definition of $u_0$, we have
  $$\nabla\log u_0(x)=-c_2 x+\nabla\log\rho(x).$$
Thus, it follows from Lemma \ref{prep-lem-5} that
  \begin{align*}
  \Big\<\nabla\log u_0(y)- \nabla\log u_0(x),\frac{y-x}{|y-x|}\Big\>
  &= -c_2|y-x|+\Big\<\nabla\log\rho(y)- \nabla\log\rho(x),\frac{y-x}{|y-x|}\Big\>\\
  &\leq 2\kappa (\log\tilde\phi_0)'\bigg(\frac{|x-y|}2\bigg).
  \end{align*}
The proof is complete.
\end{proof}

Finally, note that the function $(\log\tilde\phi_0)'$ explodes at the endpoints of the interval $[-D/2,D/2]$. For some technical reasons, we need to construct a family of smooth functions $\psi_\kappa\in C^\infty([0,D/2])$ which converges as $\kappa\to 1$ to $(\log\tilde\phi_0)'$ pointwise on the interval $[0,D/2)$.

\begin{proposition}\label{prep-prop-3}
For any $\kappa$ in a small left neighborhood of $1$, there exists a smooth function $\psi_\kappa\in C^\infty([0,D/2])$ satisfying $\psi_\kappa(0)=0$ and
\begin{itemize}
\item[\rm(i)] $\psi''_\kappa+2\psi_\kappa \psi'_\kappa-\tilde V'=0$ on $[0,D/2)$;
\item[\rm(ii)] as $\kappa\uparrow 1$, $\psi_\kappa$ converges to $(\log\tilde\phi_0)'$ pointwise on $[0,D/2)$;
\item[\rm(iii)] $\psi_\kappa$ is a modulus of log-concavity of $u_\kappa$ constructed in Proposition \ref{prep-prop-2}.
\end{itemize}
\end{proposition}

\begin{proof}
We may assume $\tilde V$ is nonnegative on $[-D/2,D/2]$, since adding a positive constant $C$ to $\tilde V$ does not change the eigen-functions $\{\tilde\phi_i\}_{i\geq 0}$, with the corresponding eigenvalues $\{\tilde\lambda_i+C\}_{i\geq 0}$. We follow the idea in the proof of \cite[Proposition 3.2]{Andrews}. Note that $\psi:= (\log\tilde\phi_0)'$ satisfies $\psi(0)=0$ and $\psi'=\tilde V-\tilde\lambda_0-\psi^2$ on the interval $(-D/2,D/2)$. Letting $q=\arctan\psi$, then $q(0)=0$ and
  $$q'-(\tilde V-\tilde\lambda_0)\cos^2 q+\sin^2 q=0\quad \mbox{on } [-D/2, D/2].$$
Here the derivatives of $q$ at the endpoints are understood as right or left derivative, respectively. For $\kappa\in (0,1)$, we consider the first order ODE
  \begin{equation}\label{prop-2.1}
  \begin{split}
  \frac{\partial q}{\partial z}-\big[(\tilde V/\kappa) -\tilde\lambda_{0}\big]\cos^2 q+\kappa \sin^2 q &=0, \quad |z|\leq D/2;\\
  \quad q(0,\kappa)&=0.
  \end{split}
  \end{equation}
By the ODE comparison theorem (cf. \cite[p.23, Theorem 8]{BR}), $q$ is strictly decreasing in $\kappa$ for $z>0$. The choice $\kappa= 1$ corresponds to $q=\arctan\psi$; hence $q(D/2,1)=-\pi/2$ and $q(z,1)\in (-\pi/2,\pi/2)$ for $0\leq z<D/2$.

Because $q(D/2,\kappa)$ is strictly decreasing in $\kappa$ and $q(D/2,1)=-\pi/2$, for $\kappa<1$ in a small neighborhood of $1$, we have $q(D/2,\kappa)> -\pi/2$ and $q(z,\kappa) \in(-\pi/2,\pi/2)$ for $0\leq z<D/2$. Set $\psi_\kappa(z)=\kappa\tan q(z,\kappa),\, |z|\leq D/2$. Direct calculations show that $\psi_\kappa$ satisfies (i). Since $\kappa<1$,
  $$\tan q(z,\kappa)>\tan q(z,1)=\psi(z) =(\log\tilde\phi_{0})'(z) \quad \mbox{for } z\in(0,D/2).$$
Thus by Proposition \ref{prep-prop-2}, $\psi_\kappa$ is a modulus of concavity of $\log u_\kappa$ for all $\kappa$ in a small left neighborhood of $1$, which means (iii) holds. It remains to check the second assertion. By the continuous dependence on $\kappa$ of the solution $q(z,\kappa)$ (cf. \cite[p.107, Corollary]{BR}), for any $z\in(0,D/2)$, $\tan q(z,\kappa)\to\tan q(z,1)$ as $\kappa$ grows to 1, which implies the pointwise convergence of $\psi_\kappa$ to $(\log\tilde\phi_{0})'$ on $[0,D/2)$. Therefore, assertion (ii) is also verified.
\end{proof}

\medskip

\noindent\textbf{Acknowledgements.} The authors are very grateful to Professors Kai He, Yong Liu, Yuan Liu and Yongsheng Song for their helpful discussions, and to Professor Elton P. Hsu for his valuable comments on the first version of this paper. They also want to thank Abraham Ng from the University of Sydney for pointing out that the function $u_0$ constructed in Proposition A.7 of the previous version has zero-gradient on the boundary $\partial\Omega$. The last two authors would like to thank Professor Liming Wu for drawing their attention to the references \cite{Carlen, MeyerZheng} and for his many suggestions on further study.

\end{document}